\documentclass[number,citesort,seceqn,dvips]{arxbj}

\aid{0}
\volume{17}
\issue{2}
\pubyear{2011}
\firstpage{814}
\lastpage{826}
\doi{10.3150/10-BEJ292}

\makeatletter
\newtheorem{theorem}{Theorem}[section]
\newtheorem{prop}[theorem]{Proposition}

\newremark{remark}[theorem]{Remark}

\newproclaim{hyp}[theorem]{Assumption}

\def\afrac#1#2{#1/(#2)}

\newcommand{\eqref}[1]{(\ref{#1})}
\newcommand{\R}{\mathbb{R}}
\renewcommand{\Re}{{\mathfrak{Re}}}

\newcommand{\Q}{\mathbb{Q}}
\newcommand{\E}{\mathbb{E}}
\newcommand{\F}{\mathcal{F}}

\newcommand{\e}{{\mathbf{e}}}
\newcommand{\ed}{\stackrel{\mathrm{(d)}}{=}}
\newcommand{\I}{\mathcal{I}}

\makeatother

\begin{document}
\begin{frontmatter}

\title{A refined factorization of the exponential law}
\runtitle{A refined factorization of the exponential law}

\begin{aug}
\author{\fnms{P.} \snm{Patie}\corref{}\ead[label=e1]{ppatie@ulb.ac.be}}

\runauthor{P. Patie}

\address{Department of Mathematics,
Universit\'e Libre de Bruxelles,
B-1050 Bruxelles,
Belgium. \\
\printead{e1}}
\end{aug}

\received{\smonth{10} \syear{2008}}
\revised{\smonth{5} \syear{2010}}

%
\begin{abstract}
Let $\xi$ be a (possibly killed) subordinator with Laplace exponent
$\phi$ and denote by $I_{\phi}=\int_0^{\infty}\mathrm{e}^{-\xi_s}\,\mathrm{d}s$, the
so-called exponential functional. Consider the positive random variable
$I_{\psi_1}$ whose law, according to Bertoin and Yor
[\textit{Electron.~Comm.~Probab.}~\textbf{6} (2001) 95--106], is determined by its
negative entire moments as follows:
\[
\E[I_{\psi_1}^{-n}] = \prod_{k=1}^n \phi(k), \qquad n=1,2,\ldots.
\]
In this note, we show that $I_{\psi_1}$ is a positive self-decomposable
random variable whenever the L\'evy measure of $\xi$ is absolutely
continuous with a monotone decreasing density. In fact, $I_{\psi_1}$ is
identified as the exponential functional of a spectrally negative
({\textit{sn}}, for short) L\'evy process. We deduce from Bertoin and
Yor [\textit{Electron.~Comm.~Probab.}~\textbf{6} (2001) 95--106] the
following factorization of the exponential law $\e$:
\[
I_{\phi} / I_{\psi_1}\ed\e,
\]
where $ I_{\psi_1}$ is taken to be independent of $I_{\phi}$. We
proceed by showing an identity in distribution between the entrance law
of an sn self-similar positive Feller process and the reciprocal of the
exponential functional of sn L\'evy processes.
As a by-product, we obtain some new examples of the law of the
exponential functionals, a new factorization of the exponential law and
some interesting distributional properties of some random variables.
For instance, we obtain that $S(\alpha)^{\alpha}$ is a
self-decomposable random variable, where $S(\alpha)$ is a positive
stable random variable of index $\alpha\in(0,1)$.
\end{abstract}

%

\end{frontmatter}
%
\section{Introduction}\label{Introduction}
Let $\xi=(\xi_t,t\geq0)$ be a possibly killed subordinator starting
from $0$, that is,~a $[0,\infty)$-valued ($\infty$ serves as absorbing
state) L\'evy process such that $\xi_0=0$. The law of $\xi$ is well
known to be characterized by its Laplace exponent $\phi$, which admits
the following L\'evy--Khintchine representation for any $u\geq0$:
%
\begin{eqnarray}\label{eq:ls}
\phi(u) = bu + \int_{0}^{\infty} (1-\mathrm{e}^{-u r})\nu(\mathrm{d}r)+q,
\end{eqnarray}
where $q\geq0$ is the killing rate, $b\geq0$ is the drift and the L\'evy measure
$\nu$ satisfies the integrability condition $\int_{\R^+} (1
\wedge r ) \nu(\mathrm{d}r) < \infty$. Note that functions of the form \eqref
{eq:ls} are also named, in the literature, as Bernstein functions. We
refer to the monographs \cite{Bertoin-96,Kyprianou-06} (resp.,
\cite{Berg-Forst-75,Jacob-01}) for a detailed account on L\'evy
processes (resp., Bernstein functions). Next, consider the so-called
exponential functional associated to $\xi$, which is defined as
\begin{eqnarray*}
I_{\phi} &=& \int_0^{\e_q} \mathrm{e}^{-\xi_s}\,\mathrm{d}s,
\end{eqnarray*}
where $\e_q$ is an independent exponential random variable with
parameter $q$ and we understand that $\e_{0}=+\infty$. Note that, for
any $q\geq0$, $I_{\phi} <\infty$ a.s.
We refer to the survey paper of Bertoin and Yor \cite{Bertoin-Yor-05}
for a thorough description of the properties of this positive random
variable and of the motivations for studying its law.
In particular, we mention that the law of $I_{\phi}$ has been
determined through its positive entire moments by Carmona \textit{et al.}~\cite{Carmona-Petit-Yor-94} as follows:
%
\begin{eqnarray}\label{eq:ms}
\E[I_{\phi}^n] &=&\frac{\Gamma(n+1)}{\prod_{k=1}^{n}\phi(k)},
\qquad n=1,2,\ldots,
\end{eqnarray}
where $\Gamma$ stands for the gamma function. Bertoin and Yor \cite{Bertoin-Yor-01} (see also \cite{Bertoin-Yor-05}, Theorem~2, for the
case $q>0$) showed that there exists a positive random variable $J$
whose law is determined by its positive entire moments as follows:
\begin{eqnarray*}
\E[J^{n}] &=&\prod_{k=1}^{n}\phi(k),\qquad n=1,2,\ldots,
\end{eqnarray*}
such that, when $J$ is taken to be independent of $I_{\phi}$, one has
the following factorization of the exponential law:
%
\begin{eqnarray} \label{eq:fe}
I_{\phi}J \stackrel{\mathrm{(d)}}{=} {\mathbf{e}},
\end{eqnarray}
where ``$\ed$'' means identity in distribution and ${\mathbf{e}}= \e_1$.

Let us now point out that since the random variable $J$ is defined on
the half-line and its law is uniquely determined by its positive entire
moments, the sequence $(s_n=\prod_{k=1}^{n}\phi(k))_{n\geq0}$
corresponds to a determinate normalized Stieltjes moment sequence. In
this direction, we should mention that Berg \cite{Berg-05} generalizes
the above fact by showing that for any $c>0$, the sequence
$(s_n^c)_{n\geq0}$ associated to a measure on the half-line $\rho_c$ is
also a Stieltjes moment sequence which is determinate for $c\leq2$. He
then deduces that there exists a unique product convolution semigroup
$(\rho_c)_{c>0}$ such that the moments of $\rho_c$ are given by $s^c_n$
for any $c>0$. Moreover, in \cite{Berg-07}, Berg characterizes the set
of normalized Stieltjes moment sequences for which this power stability
property still holds. In the same vein, Berg and Dur\'an \cite{Berg-Duran-04} study a more general mapping which allows, in
particular, the construction of a Stieltjes moment sequence of the form
$(s_n)_{n\geq0}$ with the Bernstein function $\phi$ replaced by a
completely monotone function.

The first aim of this note is to show that the random variable $1/J$ is
actually a positive self-decomposable random variable, provided that
the L\'evy measure $\nu$ in \eqref{eq:ls} admits a monotone decreasing
density. This will be achieved by identifying the random variable $1/J$
as the exponential functional of a spectrally negative L\'evy process
which we now introduce.
Let $\Xi=(\Xi_t,t\geq0)$ be a conservative spectrally negative L\'evy
process with a non-negative mean $m$ and starting from $0$, that is,~a
L\'evy process having only negative jumps such that $0\leq m=\E[\Xi
_1]<\infty$. Its law is characterized by its Laplace exponent $\psi$
which admits, in this case, the following L\'evy--Khintchine
representation for any $u\geq0$:
%
\begin{eqnarray}\label{eq:lsn}
\psi(u) = \sigma u^2 + m u + \int^{0}_{-\infty} (\mathrm{e}^{u r} -1
-ur)\Pi(\mathrm{d}r),
\end{eqnarray}
where $\sigma\geq0$ is the Gaussian coefficient and the L\'evy measure
satisfies the condition $\int_{-\infty}^0 (|r|
\wedge r^2 ) \Pi(\mathrm{d}r) <\infty$.
The exponential functional associated to $\Xi$, denoted by $I_{\psi}$,
is finite a.s.~whenever $m>0$. Its law has been determined through its
negative entire moments by Bertoin and Yor \cite{Bertoin-Yor-02}, as follows:
%
\begin{eqnarray} \label{eq:msn}
\E[I_{\psi}^{-n}] &=&m \frac{\prod_{k=1}^{n-1}\psi(k)}{\Gamma(n)},
\qquad n=1,2,\ldots,
\end{eqnarray}
with the convention that the right-hand side is $m$ when $n=1$.

We now recall that Lamperti \cite{Lamperti-72}, interested in limit
theorems for stochastic processes, shows, in particular, that for any
$x>0$, the process $X=(X_t,t\geq0)$, defined for any $t\geq0$ by
%
\begin{eqnarray} \label{eq:ss}
X_t &=& x\exp(\Xi_{A_{t/x}} ),\qquad A_t=\inf\biggl\{s\geq0; \int_0^s\mathrm{e}^{\Xi _u}\,\mathrm{d}u>t\biggr\},
\end{eqnarray}
starting from $x$ at time $0$, is a self-similar Feller process on
$(0,\infty)$ having only negative jumps. The Lamperti transformation is
actually one-to-one and extends to any L\'evy process. Bertoin and Yor
\cite{Bertoin-Yor-02}, Proposition 1, shows that the family of
probability measures $(\Q^{(\psi)}_x)_{x>0}$ of $X$, as defined in
\eqref{eq:ss}, converges as $x \downarrow0$, in the sense of
finite-dimensional distributions, to a probability measure $\Q^{(\psi
)}_0$; see also 
\cite{Caballero-Chaumont-06-b}
for the weak convergence in the Skorokhod topology. Thus, $X$ is a
Feller process on $[0,\infty)$ and Bertoin and Yor determine the law of
the random variable $J_{\psi}=(X_1,\Q^{(\psi)}_0)$, the entrance law of
$X$ at time $1$, in terms of its positive entire moments as follows:
%
\begin{eqnarray} \label{eq:me}
\E[J_{\psi}^{n}] &=& \frac{\prod_{k=1}^{n}\psi(k)}{\Gamma(n+1)},
\qquad n=1,2,\ldots.
\end{eqnarray}
They also deduce, in the case where $m>0$ and $\xi$ is the ascending
ladder height process of the dual process of $\Xi$ (see, e.g.,~\cite{Bertoin-96}, Chapter VI), that the random variable $J$, in \eqref
{eq:fe}, is $J_{\psi}$, that is,
%
\begin{eqnarray}\label{eq:feb}
I_{\phi} J_{\psi}\stackrel{\mathrm{(d)}}{=} \e.
\end{eqnarray}
The second aim of this note is to relate, in a simple way, the law of
$J_{\psi}$, for any $m\geq0$, with the exponential functional of a
spectrally negative L\'evy process.

Finally, as observed by Rivero \cite{Rivero-05}, the study of the
exponential functional is also motivated by its connection to some
interesting random equations. Indeed, from the strong Markov property
for L\'evy processes, which entails that for any finite stopping time
$T$ in the natural filtration $(\F_t,t\geq0)$ of $\xi$, the process
$(\xi_{t+T}-\xi_T, t\geq0)$ is independent of $\F_T$ and has the same
distribution as $\xi$, we readily deduce that the random variable
$I_{\phi}$, in the case $q=0$, is a solution to the random affine equation
%
\begin{eqnarray} \label{eq:re}
I_{\phi} \stackrel{\mathrm{(d)}}{=} \int_0^{T} \mathrm{e}^{-\xi_s}\,\mathrm{d}s + \mathrm{e}^{-\xi_T}I'_{\phi},
\end{eqnarray}
where, on the right-hand side, $I'_{\phi}$ is an independent copy of
$I_{\phi}$. Note that this type of random equation have been studied by
Kesten \cite{Kesten-73} and Goldie \cite{Goldie-91}. By means of a
similar argument, but using the absence of positive jumps of $\Xi$ (see
\cite{Rivero-03}, Proposition 4, for more details), we get that
$I_{\psi}$ is a solution to the random affine equation, for any $y>0$,
%
\begin{eqnarray} \label{eq:rec}
I_{\psi} \stackrel{\mathrm{(d)}}{=} \int_0^{T_{y}} \mathrm{e}^{-\Xi_s}\,\mathrm{d}s + \mathrm{e}^{-y}I'_{\psi},
\end{eqnarray}
where $T_{y}=\inf\{s>0; \Xi_s\geq y\}$ and, on the right-hand side,
$I'_{\psi}$ is an independent copy of $I_{\psi}$. Hence, $I_{\psi}$ is
a positive self-decomposable random variable and, in particular, its
law is absolutely continuous and unimodal; see, e.g., 
\cite{Sato-99} and
\cite{Steutel-vanHarn-04} for an
excellent account of this set of probability measures.

\section{Main results}\label{sec2}

\subsection{Factorization of the exponential law with exponential functionals}
In this subsection, we suppose that $\xi$ is a subordinator starting
from $0$ with Laplace exponent given by \eqref{eq:ls}. We introduce the
following hypothesis on the L\'evy measure of~$\xi$.
\begin{hyp} \label{hyp1}
There exists a monotone decreasing function $f$ such that $\nu
(\mathrm{d}x)=f(x)\,\mathrm{d}x.$
\end{hyp}

We recall that, under this condition, $-\mathrm{d}f(x)$ is a Stieltjes measure
on $(0,\infty)$. We also use the notation $-\mathrm{d}f(-x)$ for the image of
the positive measure $-\mathrm{d}f(x)$ under the map $x\mapsto-x$. For instance,
if $f$ is, in addition, differentiable, then $-\mathrm{d}f(-x)=-f'(-x)\,\mathrm{d}x$.
We are now ready to derive our refinement of the factorization of the
exponential law.
\begin{theorem} \label{thm:1}
Let $\xi$ be a subordinator with Laplace exponent $\phi$ given by \eqref
{eq:ls}. If Assumption~\ref{hyp1} holds, then there exists an
independent spectrally negative L\'evy process with a positive mean and
Laplace exponent $\psi_1$, analytic in the domain $C=\{u \in\mathbb{C}; \Re(u)> -1\}$, with $\psi_1(-1)=-\phi(0)$, given by
\begin{eqnarray*}
\psi_1(u)=bu^2+\phi(1)u +\int^{0}_{-\infty}(\mathrm{e}^{ur}-1-ur)\Pi(\mathrm{d}r),\qquad u \in C,
\end{eqnarray*}
where $\Pi(\mathrm{d}r)=\mathrm{e}^r (f(-r)\,\mathrm{d}r-\mathrm{d}f(-r) )$ is a Stieltjes measure on
$(-\infty,0)$. Moreover, the law of the positive self-decomposable
random variable $I_{\psi_1}$ is determined by its negative entire
moments as follows:
\begin{eqnarray*}
\E[I_{\psi_1}^{-n}] &=& \prod_{k=1}^{n}\phi(k), \qquad n=1,2,\ldots.
\end{eqnarray*}
The exponential law admits the following factorization:
%
\begin{eqnarray}\label{eq:nfe}
I_{\phi}/I_{\psi_1}
\stackrel{\mathrm{(d)}}{=} {\mathbf{e}},
\end{eqnarray}
where ${\mathbf{e}}$ stands for an exponential random variable of
parameter $1$.

Conversely, if $\psi$ is of the form \eqref{eq:lsn} with $m>0$ and is
analytic in the domain $C$ with $\psi(-1)\leq0$, then there exists an
independent subordinator with Laplace exponent $\phi_{-1}$ given by
\begin{eqnarray*}
\phi_{-1}(u)=-\psi(-1)+\sigma u +\int^{\infty}_{0}(1-\mathrm{e}^{-ur})\mathrm{e}^{r}\Pi
(-\infty,-r)\,\mathrm{d}r,\qquad u\geq0,
\end{eqnarray*}
such that
%
\begin{eqnarray}\label{eq:nfe2}
I_{\phi_{-1}}/I_{\psi}
\stackrel{\mathrm{(d)}}{=} {\mathbf{e}}.
\end{eqnarray}

\end{theorem}
\begin{remark}
(1) We have several comments to offer on the identity \eqref{eq:nfe}
when compared to \eqref{eq:feb}.
First, our hypotheses are slightly less restrictive. Indeed, it is well
known (see, e.g.,~\cite{Bertoin-96}, Chapter~VI) that the ascending
ladder height process of the dual process of $\Xi$ satisfies Assumption
\ref{hyp1} and is a killed subordinator; thus, in \eqref{eq:feb}, $q$
is necessarily positive. More importantly, we have identified the
mixture random variable of $I_{\phi}$ in the factorization of the
exponential law as the reciprocal of a positive self-decomposable
random variable. Finally, our identity allows further explicit examples
to be obtained for the law of the exponential functional of L\'evy
processes. All of these facts will be illustrated in Section~\ref{sec3}.
\smallskipamount=0pt
\begin{longlist}[(1)]
\item[(2)] The analyticity property of $\psi_1$ means that the associated
spectrally negative L\'evy process $\Xi^{1}$ admits exponential moments
of order $u\geq-1$, that is,~ for any $u\geq-1$, we have
\[
\E[\mathrm{e}^{u\Xi^{1}_1}]<\infty.
\]
We shall show, in Proposition~\ref{prop:1} below, how to construct a
spectrally negative L\'evy process with such a property from any
spectrally negative L\'evy process with a non-negative mean.
\end{longlist}
\end{remark}
\begin{pf*}{Proof of Theorem~\ref{thm:1}}
Let us write $\psi_1(u)=u\phi(u+1)$. Then, recalling that $\phi$ is
analytic in the right half-plane, we readily deduce that the mapping
$u\mapsto\psi_1(u)$ is analytic in $C$, with $\psi_1(0)=0$ and $\psi
_1(-1)=-\phi(0)$. Let us now show that under Assumption~\ref{hyp1},
$\psi_1$ is the Laplace exponent of a spectrally negative L\'evy
process with a positive mean. On one hand, since $r\mapsto f(r)$ is
monotone decreasing on $(0,\infty),$ it follows that $\Pi(\mathrm{d}r)=\mathrm{e}^r
(f(-r)\,\mathrm{d}r-\mathrm{d}f(-r) )$ is clearly a Stieltjes measure on $(-\infty,0)$.\vadjust{\goodbreak}
On the other hand, by integration by parts and a change of variable, we
have, for any $u\geq0$,
\begin{eqnarray*}
\psi_1(u)
&=& u \biggl(b(u+1)+ \int_{0}^{\infty} \bigl(1-\mathrm{e}^{-(u+1)r}\bigr)f(r)\,\mathrm{d}r+q \biggr)\\
&=& bu^2 + \biggl(b+q+ \int_{0}^{\infty}(1-\mathrm{e}^{-r})f(r)\,\mathrm{d}r \biggr) u + u\int_0^{\infty}(1-\mathrm{e}^{-ur})\mathrm{e}^{-r}f(r)\,\mathrm{d}r\\
&=& bu^2 +\phi(1)u + \int_0^{\infty}(\mathrm{e}^{-ur}-1+ur)\mathrm{e}^{-r}\bigl(f(r)\,\mathrm{d}r-\mathrm{d}f(r)\bigr)\\
&=&bu^2+\phi(1)u +\int^{0}_{-\infty}(\mathrm{e}^{ur}-1-ur)\Pi(\mathrm{d}r).
\end{eqnarray*}
Checking, by integration by parts, that $\int^0_{-\infty} (|r|\wedge
r^2) \Pi(\mathrm{d}r)<\infty$, we get that $\psi_1$ is the Laplace exponent of a
spectrally negative L\'evy process with a positive mean since $\psi
_1(0^+)=\phi(1)>0$. Then, by means of the identity \eqref{eq:msn}, we
obtain, for any $n=1,2,\ldots,$
\begin{eqnarray*}
\E[I_{\psi_1}^{-n}] =\phi(1) \frac{\prod_{k=1}^{n-1}\psi_1(k)}{\Gamma
(n)}
= \phi(1)\frac{\prod_{k=1}^{n-1}k\phi(k+1)}{\Gamma(n)}
= \prod_{k=1}^{n}\phi(k),
\end{eqnarray*}
where we have used the identity $\Gamma(n)=\prod_{k=1}^{n-1}k$. The
self-decomposability of $I_{\psi_1}$ was discussed in the \hyperref[Introduction]{Introduction}
and the factorization of the exponential law follows readily from the
independence of the random variables $I_{\psi_1}$ and $I_{\phi}$ and
the identity \eqref{eq:ms}. The converse follows by means of similar
reasoning. We only need to check that $\phi_{-1}(u)=\frac{\psi
(u-1)}{u-1}$ is the Laplace exponent of a subordinator. From the
conditions imposed on $\psi$, we can easily deduce that $\phi_{-1}$ is
well defined on $\R^+$ with $\phi_{-1}(0)=-\psi(-1)\geq0$. Moreover, by
integrations by parts,
we get
\begin{eqnarray*}
\phi_{-1}(u) &=& \sigma u -(\sigma-m)+ \frac{1}{u-1}\int_{-\infty}^0
\bigl(\mathrm{e}^{(u-1)r}-1-(u-1)r\bigr)\Pi(\mathrm{d}r)\\
&=& \sigma u -(\sigma-m)- \int_{-\infty}^0 \bigl(\mathrm{e}^{(u-1)r}-1\bigr)\Pi(-\infty,r)\,\mathrm{d}r\\
&=& \sigma u -(\sigma-m) - \int^{0}_{-\infty} (\mathrm{e}^{ur}-1)\mathrm{e}^{-r}\Pi
(-\infty,r)\,\mathrm{d}r-\int^{0}_{-\infty} (\mathrm{e}^{-r}-1+r)\Pi(\mathrm{d}r)\\
&=& -\psi(-1)+\sigma u +\int_{0}^{\infty} (1-\mathrm{e}^{-ur})\mathrm{e}^{r}\Pi(-\infty,-r)\,\mathrm{d}r.
\end{eqnarray*}
The proof of the theorem is thus completed.
\end{pf*}

As a direct consequence of Theorem~\ref{thm:1}, we have the following
fact. By self-decomposability, the law of $I_{\psi_1}$ is absolutely
continuous and thus the random variable $J$ in \eqref{eq:fe} admits a
density with respect to the Lebesgue measure, which, according to
Bertoin and Yor \cite{Bertoin-Yor-01}, is a $1$-harmonic function for
the self-similar process associated to $\xi$ in the Lamperti mapping
\eqref{eq:ss}. Thus, writing $p_{1}$ for the density of $I_{\psi_1}$,
the mapping $x\mapsto x^{-2}p_1 (x^{-1} )$ on $\R^+$ is a $1$-harmonic
function for the self-similar process associated to $\xi$.

We complete this part with the following observation. Let us suppose
that there exists a subordinator with Laplace exponent $\widehat{\phi}$
such that
\[
\phi(u)\widehat{\phi}(u)=u, \qquad u\geq0.
\]
Under such a condition, $\phi$ is called a special Bernstein function
and we refer to Kyprianou and Rivero \cite{Kyprianou-Rivero-08}, and
references therein, for more information on this function. Note that
such an identity occurs in fluctuation theory for L\'evy processes and
that sufficient conditions for $\phi$ to be a special Bernstein
function are that Assumption~\ref{hyp1} holds, $d=0$ and the mapping
$x\mapsto f(x) (\int_x^{\infty}f(y)\,\mathrm{d}y )^{-1}$ is a decreasing function
on $(0,\infty)$.
As remarked in \cite{Bertoin-Yor-01}, we have, in this case, with the
obvious notation,
\begin{eqnarray*}
I_{\widehat{\phi}}I_{\phi}
\stackrel{\mathrm{(d)}}{=} \e.
\end{eqnarray*}
We deduce readily, from Theorem~\ref{thm:1}, that if $\phi$ is a
special Bernstein function satisfying Assumption~\ref{hyp1}, then
\begin{eqnarray*}
I_{\widehat{\phi}}
\stackrel{\mathrm{(d)}}{=} 1/I_{\psi_1}.
\end{eqnarray*}
Moreover, if, in addition, $\widehat{\phi}$ also satisfies Assumption
\ref{hyp1}, then $ 1/I_{\phi}$ is a positive self-decomposable random variable.
Note that if $\widehat{\phi}(0)=0$, then $I_{\widehat{\phi}}$ is the
solution to the random affine equation \eqref{eq:re} and thus, in this
case, solving this equation reduces to solve the random affine equation
with constant coefficient \eqref{eq:rec}.
\subsection{Exponential functionals and entrance laws $J_{\psi}$}

We now assume that $\Xi$ is a spectrally negative L\'evy process with a
non-negative mean $m$. Its Laplace exponent has the form \eqref{eq:lsn}.
We recall that the positive random variable $J_{\psi}$ is the entrance
law at time $1$ of the self-similar Feller process associated to $\Xi$
via the Lamperti mapping \eqref{eq:ss}. We also mention that, when
$m>0$, Bertoin and Yor \cite{Bertoin-Yor-02-b} show that the
distribution of $1/I_{\psi}$ is the so-called length-biased
distribution of $J_{\psi}$, that is,~using the identity $\E[1/I_{\psi}]=m$,
\[
\E[g(J_{\psi})]= m^{-1}\E[1/I_{\psi}g(1/I_{\psi})]
\]
for any measurable function $g\dvtx \R^+ \rightarrow\R^+$. We refine this
connection in the following proposition.
\begin{prop} \label{prop:1}
Let $\psi$ be the Laplace exponent of a spectrally negative L\'evy
process with mean $m\geq0$. Then, the mapping defined by
\[
\psi_2(u)=\frac{u}{u+1}\psi(u+1)
\]
is analytic in $C=\{u\in\mathbb{C}; \Re(u)>-1\}$ and is the Laplace
exponent of a spectrally negative L\'evy process with a positive mean
$\psi(1)$. Moreover, the identity in distribution
\begin{eqnarray*}
J_{\psi}
\stackrel{\mathrm{(d)}}{=} 1/I_{\psi_2}
\end{eqnarray*}
holds. Consequently, the law of $J_{\psi}$ is absolutely continuous for
any $m\geq0$.
\end{prop}
\begin{pf}
First, since it is well known that $\psi$ is analytic in the right
half-plane, it is clear that the mapping $\psi_2(u)=\frac{u}{u+1}\psi
(u+1)$ is analytic in $C$. Next, we recall that $\psi$ has the form
\begin{eqnarray*}
\psi(u)= \sigma u^2 + mu+ \int_{-\infty}^{0}(\mathrm{e}^{ur}-1-ur)\Pi(\mathrm{d}r),
\end{eqnarray*}
where $\int^0_{-\infty}(|r|\wedge r^2)\Pi(\mathrm{d}r)<\infty$ and $\sigma\geq
0$. Thus, by means of integration by parts and writing $f(r)=\Pi(-\infty
,r)$ for the tail of the L\'evy measure, we get
%
\begin{eqnarray} \label{eq:p1}
&&\frac{u}{u+1}\psi(u+1)\nonumber
\\
&&\quad= \sigma u^2 + (m+\sigma)u+ \frac{u}{u+1}\int_{-\infty}^0\bigl(\mathrm{e}^{(u+1)r}-1-(u+1)r\bigr)\Pi(\mathrm{d}r)\nonumber\\
&&\quad= \sigma u^2 + (m+\sigma)u- u\int_{-\infty}^0\bigl(\mathrm{e}^{(u+1)r}-1\bigr)f(r)\,\mathrm{d}r\nonumber\\
&&\quad= \sigma u^2 + (m+\sigma)u- u \biggl(\int_{-\infty}^{0}(\mathrm{e}^{ur}-1)\mathrm{e}^rf(r)\,\mathrm{d}r
+\int_{-\infty}^0(\mathrm{e}^{r}-1)f(r)\,\mathrm{d}r \biggr)
\\
&&\quad= \sigma u^2 + \biggl(m+\sigma+\int_{-\infty}^0(\mathrm{e}^{r}-1-r)\Pi(\mathrm{d}r)\biggr)u\nonumber\\
&&\qquad{}+\int_{-\infty}^0(\mathrm{e}^{ur}-1-ur)\mathrm{e}^r\bigl(f(r)\,\mathrm{d}r+\Pi(\mathrm{d}r)\bigr)\nonumber\\
&&\quad= \sigma u^2 + \psi(1)u+\int_{-\infty}^{0}(\mathrm{e}^{ur}-1-ur)\mathrm{e}^r\bigl(f(r)\,\mathrm{d}r+\Pi(\mathrm{d}r)\bigr),\nonumber
\end{eqnarray}
where we recognize the Laplace exponent of a spectrally negative L\'evy
process. Finally, observing that $\lim_{u\rightarrow0}\frac{\mathrm{d}}{\mathrm{d}u}\frac
{u}{u+1}\psi(u+1)=\psi(1)>0$ since $\psi$ is increasing on $(0,\infty
)$, we have, from \eqref{eq:msn} and any $ n=1,2\ldots,$
\begin{eqnarray*}
\E[I_{\psi_2}^{-n}]
=\psi(1)\frac{\prod_{k=1}^{n-1}(\afrac{k}{k+1})\psi(k+1)}{\Gamma(n)}
=\frac{\prod_{k=1}^{n}\psi(k)}{\Gamma(n+1)}
=\E[J_\psi^{n}],
\end{eqnarray*}
where the last identity follows from \eqref{eq:me}. The absolute
continuity property of the law of $J_{\psi}$ follows from that of
$I_{\psi_2}$ as a self-decomposable random variable.
\end{pf}
We mention that the random variable $J_{\psi}$ appears in the study of
the so-called Ornstein--Uhlenbeck process associated to $X$. Indeed, if
one considers the stochastic process $U=(U_t,t\geq0)$ defined, for any
$t\geq0$, by
\[
U_t = \mathrm{e}^{-t}X_{\mathrm{e}^{t}-1},
\]
then $U$ is a stationary Feller process on $[0,\infty)$ and its unique
invariant measure is the law of $J_{\psi}$; see, for instance,~\cite{Patie-08a}, Theorem~1.2. The above proposition tells us that the
invariant measure is absolutely continuous for any $m\geq0$.

We also indicate that the transformation of the Laplace exponent of a
spectrally negative L\'evy process used in the proof of Proposition~\ref{prop:1} is a specific instance of more general mappings of
characteristic exponents of L\'evy processes introduced and studied by
Kyprianou and Patie \cite{Kyprianou-Patie-08}.

\section{Some examples}\label{sec3}
In this section, we will make use of the identities presented in
Section~\ref{sec2} to obtain new explicit examples of the law of the exponential
functional associated to subordinators or spectrally negative L\'evy
processes, to obtain a new factorization of the exponential law and to
prove the self-decomposability property of some positive random variables.

In \cite{Bertoin-Yor-01}, the authors study the connection between the
law of the exponential functional of some subordinators and the
following factorization of the exponential law:
%
\begin{equation}\label{eq:fee}
\mathbf{e}^{\alpha}S(\alpha)^{-\alpha}\ed\e,
\end{equation}
where $\alpha\in(0,1)$ and $S(\alpha)$ is a positive $\alpha$-stable
random variable, independent of $\e$. We split this example into two parts.

On the one hand, the authors show that
\[
I_{\phi}\ed S(\alpha)^{-\alpha}
\]
with
%
\begin{eqnarray}
\phi(u)= \frac{\alpha\Gamma(\alpha u+1)}{\Gamma(\alpha(u-1)+1)} =\int
_0^{\infty}(1-\mathrm{e}^{-ur})f(r)\,\mathrm{d}r, \label{eq:lpp}
\end{eqnarray}
where
%
\begin{eqnarray} \label{eq:df}
f(r) =\frac{\mathrm{e}^{-r/\alpha}}{\Gamma(1-\alpha)(1-\mathrm{e}^{-r/\alpha})^{\alpha
+1}},\qquad r>0.
\end{eqnarray}
\begin{enumerate}[(1)]
\item[(1)] We start by applying the first part of Theorem~\ref{thm:1}.
It is easy to check that the mapping $r\mapsto f(r)$ is decreasing on
$(0,\infty)$ and thus, from Theorem~\ref{thm:1} and using the
recurrence relation $\Gamma(u+1)=u\Gamma(u), u>0$, we get that
\begin{eqnarray*}
\psi_1(u)&=& \frac{\Gamma(\alpha(u+1)+1)}{\Gamma(\alpha u)},
\end{eqnarray*}
which, after some easy calculations, yields
%
\begin{eqnarray} \label{eq:lm}
\psi_1(u)&=& \alpha\Gamma(\alpha+1)u+\int_{-\infty}^0(\mathrm{e}^{ur}-ur-1)\frac
{(\alpha+1)\mathrm{e}^{(\alpha+1)
r/\alpha}}{\alpha\Gamma(1-\alpha) (1-\mathrm{e}^{r/\alpha})^{\alpha+2}}\,\mathrm{d}r.
\end{eqnarray}
Hence, from \eqref{eq:fee}, we deduce that
\[
I_{\psi_1} \ed\mathbf{e}^{-\alpha}.
\]
This result, up to a multiplicative constant, was actually obtained by
Patie in \cite{Patie-CBI-09}, Theorem~4.1, where it is shown that the
law of $I_{\psi_1}$ is related to the distribution of the absorption
time of the $\alpha$-self-similar continuous-state branching process.
\item[(2)] Moreover, let us define, as in Proposition~\ref{prop:1},
\begin{eqnarray*}
\psi_2(u)= \frac{u}{u+1}\psi_1(u+1)
= \alpha u \frac{\Gamma(\alpha(u+2)+1)}{\Gamma(\alpha(u+1)+1)}
\end{eqnarray*}
and, after observing that
\[
\frac{\partial}{\partial r} \biggl(\frac{\mathrm{e}^{(\alpha+1)r/\alpha
}}{(1-\mathrm{e}^{r/\alpha})^{\alpha+1}} \biggr)= \frac{(\alpha+1)\mathrm{e}^{(\alpha+1)
r/\alpha}}{\alpha(1-\mathrm{e}^{r/\alpha})^{\alpha+2}},
\]
we obtain, from \eqref{eq:p1} and \eqref{eq:lm},
\begin{eqnarray*}
\psi_2(u)&=& \frac{\Gamma(2\alpha+1)}{\Gamma(\alpha)}u+\int^0_{-\infty
}(\mathrm{e}^{ur}-ur-1)\frac{ \mathrm{e}^{(2\alpha+1)
r/\alpha}}{\Gamma(1-\alpha) (1-\mathrm{e}^{r/\alpha})^{\alpha+2}} \biggl(\frac{2\alpha
+1}{\alpha}-\mathrm{e}^{r/\alpha} \biggr)\,\mathrm{d}r.
\end{eqnarray*}
Hence, from \eqref{eq:msn}, we get, for any $n=1,2,\ldots,$
\begin{eqnarray*}
\E[I_{\psi_2}^{-n}]= \alpha^{n-1}\frac{\Gamma(2\alpha+1)}{\Gamma
(\alpha)}\frac{\Gamma(\alpha n+\alpha+1)}{\Gamma(2\alpha+1)}
= \alpha^{n}\frac{\Gamma(\alpha n+\alpha+1)}{\Gamma(\alpha+1)}
\end{eqnarray*}
and, by moment identification, we have
\[
\alpha I_{\psi_2} \ed G^{-\alpha}(\alpha+1),
\]
where $G(a)$ stands for a gamma random variable of parameter $a>0$. We
deduce that, for any $\alpha\in(0,1)$, the random variable $G^{-\alpha
}(\alpha+1)$ is a positive self-decomposable random variable.
\item[(3)]
Next, we apply the converse part of Theorem~\ref{thm:1} to $\psi_2$. To
this end, we introduce the subordinator with Laplace exponent $\phi
_{-1}$ defined by
\begin{eqnarray*}
\phi_{-1}(u)= \frac{1}{u-1}\psi_2(u-1)
= \frac{\psi_1(u)}{u}
= \alpha \frac{\Gamma(\alpha(u+1)+1)}{\Gamma(\alpha u+1)}.
\end{eqnarray*}
To get the L\'evy--Khintchine representation of $\phi_{-1}$, note, from
\eqref{eq:lpp}, that
\[
\phi_{-1}(u)=\phi(u+1)-\phi(1)+\phi(1).
\]
That is, $\phi_{-1}$ is the Laplace exponent of the Esscher transform
of $\phi$ killed at rate $\phi(1)$. Thus, using the expression \eqref{eq:df},
\begin{eqnarray*}
\phi_{-1}(u)&=& \alpha\Gamma(\alpha+1)+ \int_0^{\infty}(1-\mathrm{e}^{-ur})\mathrm{e}^{-r}f(r)\,\mathrm{d}r.
\end{eqnarray*}
We have, from \eqref{eq:ms},
\[
\E[I_{\phi_{-1}}^{n}]= \alpha^{-n}\frac{\Gamma(\alpha+1)\Gamma
(n+1)}{\Gamma(\alpha n+ \alpha+1)},\qquad n=1,2,\ldots.
\]
In order to characterize the law of $I_{\phi_{-1}}$, let us denote by
$U$ an uniform random variable on $(0,1)$ and by $S^{-\alpha}_1(\alpha
)$ a random variable distributed according to the length-biased
distribution of $S^{-\alpha}(\alpha)$, that is,~for any measurable
function $g\dvtx\R^+\rightarrow\R^+$, we have
\[
\E[g (S_1^{-\alpha}(\alpha) ) ]= \frac{\E[S^{-\alpha}(\alpha)g
(S^{-\alpha}(\alpha) ) ]}{\E[S^{-\alpha}(\alpha) ]}.
\]
Recalling that for any $n=1,2,\ldots,$
\[
\E[S^{-\alpha n}(\alpha) ]= \frac{\Gamma(n+1)}{\Gamma(\alpha n+1)},
\]
we get, by taking the random variable $U$ independent of $S^{-\alpha
}_1(\alpha)$,
\begin{eqnarray*}
\E[U^{n}S_1^{-\alpha n}(\alpha) ]&=& \frac{\E[U^{n} ]\E[S^{-\alpha
(n+1)}(\alpha) ]}{\E[S^{-\alpha}(\alpha) ]}\\
&=& \frac{\Gamma(\alpha+1)\Gamma(n+2)}{(n+1)\Gamma(\alpha n+ \alpha
+1)}\\
&=& \frac{\Gamma(\alpha+1)\Gamma(n+1)}{\Gamma(\alpha n+ \alpha+1)}.
\end{eqnarray*}
By moment identification, we deduce the following identity:
\[
\alpha I_{\phi_{-1}} \ed US^{-\alpha}_1(\alpha).
\]
This yields the following factorization of the exponential law:
\[
US^{-\alpha}_1(\alpha) G^{\alpha}(\alpha+1)\ed\e,
\]
where the three random variables on the left-hand side are assumed to
be independent.
\end{enumerate}

On the other hand, Bertoin and Yor \cite{Bertoin-Yor-01}
also observed that
\[
I_{\widehat{\phi}}\ed\mathbf{e}^{\alpha}
\]
with
\begin{eqnarray*}
\widehat{\phi}(u)&=& u\frac{\Gamma(\alpha(u-1)+1)}{\Gamma(\alpha u+1)}
=\int_0^{\infty}(1-\mathrm{e}^{-ur})\frac{(1-\alpha)^2\mathrm{e}^{r/\alpha}}{\alpha
\Gamma(\alpha+1)(\mathrm{e}^{r/\alpha}-1)^{2-\alpha}}\,\mathrm{d}r.
\end{eqnarray*}
Easily verifying that the density is decreasing, we obtain, appealing
to obvious notation, that
\begin{eqnarray*}
\widehat{\psi}_1(u)&=& u\frac{\Gamma(\alpha u+1)}{\alpha\Gamma(\alpha(u+1))},
\end{eqnarray*}
which, after some easy manipulations, yields
\begin{eqnarray*}
\widehat{\psi}_1(u)&=& \Gamma^{-1}(\alpha+1)u
\\
&&{}+\int_{-\infty
}^0(\mathrm{e}^{ur}-ur-1)\frac{(1-\alpha)^2\mathrm{e}^{-r/\alpha
}}{\alpha^2 \Gamma(\alpha+1)(\mathrm{e}^{-r/\alpha}-1)^{2-\alpha}} \biggl(1-\alpha
+\frac{2-\alpha}{(1-\mathrm{e}^{r/\alpha})} \biggr) \,\mathrm{d}r.
\end{eqnarray*}
Thus, from the identity \eqref{eq:fee} and Theorem~\ref{thm:1}, we
deduce that
\[
\I_{\widehat{\psi}_1} \ed S(\alpha)^{\alpha}.
\]
Hence, $S(\alpha)^{\alpha}$ is a positive self-decomposable random variable.

\section*{Acknowledgements}
The author is grateful to an anonymous referee for a careful reading of this paper.

\printhistory

\end{document}